\documentclass[12pt]{article}
\usepackage{amssymb}
\begin{document}
\baselineskip=19pt
\parskip=7.5pt
\newcommand{\la}{\langle}
\newcommand{\ra}{\rangle}
\newcommand{\psp}{\vspace{0.4cm}}
\newcommand{\pse}{\vspace{0.2cm}}
\newcommand{\ptl}{\partial}
\newcommand{\dlt}{\delta}
\newcommand{\sgm}{\sigma}
\newcommand{\al}{\alpha}
\newcommand{\be}{\beta}
\newcommand{\G}{\Gamma}
\newcommand{\gm}{\gamma}
\newcommand{\vs}{\varsigma}
\newcommand{\lmd}{\lambda}
\newcommand{\td}{\tilde}
\newcommand{\vf}{\pi}
\newcommand{\rd}{\mbox{Rad}}
\newcommand{\ad}{\mbox{ad}}
\newcommand{\stl}{\stackrel}
\newcommand{\ul}{\underline}
\newcommand{\es}{\epsilon}
\newcommand{\dmd}{\diamond}
\newcommand{\clt}{\clubsuit}
\newcommand{\vt}{\vartheta}
\newcommand{\ves}{\varepsilon}
\newcommand{\dg}{\dagger}
\newcommand{\kn}{\mbox{ker}}
\newcommand{\im}{\mbox{im}}
\newcommand{\for}{\mbox{for}}
\newcommand{\dvg}{\mbox{div}}
\newcommand{\rar}{\rightarrow}
\newcommand{\NJ}{\mathbb{N}^{\ell_1+\ell_2}}
\newcommand{\ZJ}{\mathbb{Z}^{\ell_1+\ell_2}}
\newcommand{\bs}{\backslash}
\newcommand{\der}{\mbox{Der}\:}
\newcommand{\lra}{\Longleftrightarrow}
\newcommand{\BB}{\mathbb}
\def\a{\alpha}
\def\b{\beta}
\def\th{\theta}
\def\N{\hbox{$I\hskip -4pt N$}}
\def\F{\hbox{$I\hskip -4pt F$}}
\def\Z{\hbox{$Z\hskip -5.2pt Z$}}
\def\sZ{\hbox{$\sc Z\hskip -4.2pt Z$}}
\def\Q{\hbox{$Q\hskip -5pt\vrule height 6pt depth 0pt\hskip 6pt$}}
\def\sQ{\hbox{$\sc Q\hskip -4pt\vrule height 6pt depth 0pt\hskip 6pt$}}
\def\R{\hbox{$I\hskip -3pt R$}}
\def\C{\hbox{$C\hskip -5pt \vrule height 6pt depth 0pt \hskip 6pt$}}
\def\N{\mathbb{N}}
\def\Z{\mathbb{Z}}
\def\sZ{\mathbb{Z}}
\def\Q{\mathbb{Q}}
\def\sQ{\mathbb{Q}}
\def\R{\mathbb{R}}
\def\sR{\mathbb{R}}
\def\C{\mathbb{C}}
\def\sC{\mathbb{C}}
\def\F{\mathbb{F}}
\def\J {\rb{5pt}{\mbox{$^{\rar}_{\dis J}$}}}
\def\qed{\hfill \hfill \ifhmode\unskip\nobreak\fi\ifmmode\ifinner
\else\hskip5pt\fi\fi\hbox{\hskip5pt\vrule width4pt height6pt
depth1.5pt\hskip 1 pt}}
\def\d{\delta}
\def\D{\Delta}
\def\g{\gamma}
\def\G{\Gamma}
\def\l{\lambda}
\def\L{\Lambda}
\def\o{\omiga}
\def\p{\psi}
\def\Si{\Sigma}
\def\si{\sigma}
\def\sc{\scriptstyle}
\def\ssc{\scriptscriptstyle}
\def\dis{\displaystyle}
\def\cl{\centerline}
\def\ll{\leftline}
\def\rl{\rightline}
\def\sF{\hbox{$\sc I\hskip -2.3pt F$}}
\def\nl{\newline}
\def\ol{\bar}
\def\ul{\underline}
\def\wt{\widetilde}
\def\wh{\widehat}
\def\rar{\rightarrow}
\def\Rar{\Rightarrow}
\def\lar{\leftarrow}
\def\Lar{\Leftarrow}
\def\rla{\leftrightarrow}
\def\Rla{\Leftrightarrow}
\def\bs{\backslash}
\def\hs{\hspace*}
\def\vs{\vspace*}
\def\rb{\raisebox}
\def\ra{\rangle}
\def\la{\langle}
\def\SS{\hbox{$S\hskip -6.2pt S$}}
\def\ni{\noindent}
\def\hi{\hangindent}
\def\ha{\hangafter}
\def\Box{\mbox{$|\!\!\!\!|\!\!\!\!|$}}
\def\AA{{\cal A}}
\def\BB{{\cal B}}
\def\CC{{\cal C}}
\def\DD{{\cal D}}
\def\ii{{\vec i}}
\def\jj{{\vec j}}
\def\kk{{\vec k}}
\def\mm{{\vec m}}
\def\HOM{\mbox{${\rm Hom}^*_{\sZ}(\BB,\F)$}}
\def\sone{{1\hskip -5.5pt 1}}
\def\one{{1\hskip -6.5pt 1}}
\def\lsb{\mbox{{\small(}}}
\def\rsb{\mbox{{\small)}}}
\def\sF{\mathbb{F}}
%
%
\begin{center}
{\large\bf Structure of a new class of non-graded infinite
dimensional simple Lie algebras}\footnote{AMS Subject
Classification-Primary: 17B40, 17B65, 17B70\nl\hs{4ex}Supported by
a NSF grant 10171064 of China and two research grants from
Ministry of Education of China}
\end{center}
\vs{-7pt}
\cl{(appeared in {\em J.~Algebra} {\bf267} (2003) 542-558)}
\vs{-6pt}\par
\cl{Yucai Su}
\par{\small\it
\cl{Department of Mathematics, Shanghai Jiaotong University,
Shanghai 200030, P.~R.~China} \cl{Email: ycsu@sjtu.edu.cn}}
\par
\ \vs{-15pt}\par
\begin{abstract}{
A new class of infinite dimensional simple Lie algebras
over a field with characteristic 0 are constructed.
These are examples
of non-graded Lie algebras. The isomorphism classes of these Lie
algebras are determined. The structure space of these algebras is given
explicitly.
}\end{abstract}
\par\
\vs{-15pt}\par
\cl{\bf\S 1. \ Introduction}
\par
Block [1] introduced a class of infinite dimensional simple Lie algebras
over a field with characteristic zero. Since then (and in particular
recently), infinite dimensional simple Lie algebras with some features of
Block algebras have received some authors' interests (see for example
[2-7,10-11,13-16,18-20]). This is probably partially because the Block
algebras are closely related to the (centerless) Virasoro algebra
(for example, [19,20]).
\par
One observes that the classical Block algebras $\BB$ are graded Lie algebras,
i.e., $\BB$ can be finitely graded: $\BB=\oplus_{\a\in\G}\BB_\a,
[\BB_\a,\BB_\b]\subset\BB_{\a+\b}$ for $\a,\b\in\G$, where $\G$ is some
Abelian group $\G$ such that all homogeneous spaces $\BB_\a$ are finite
dimensional. It is pointed out in [17] that
non-graded infinite dimensional Lie algebras appear
naturally in the theory of Hamiltonian operators, the theory of vertex
algebras and their multi-variable analogues; that they play important roles
in mathematical physics; that the structure spaces of some families of
non-graded simple Lie algebras can be viewed as analogues of vector bundles
with Lie algebras as fibers and it is also believed that some of non-graded
simple Lie algebras could be related to noncommutative geometry.
\par
The key constructional ingredients of the non-graded infinite dimensional
simple Lie algebras in [17] are locally finite derivations.
In this paper,
based on the classification [9] of the pairs $(\AA,\DD)$, where $\AA$ is a
commutative associative algebra with an identity element and $\DD$ is
a finite dimensional commutative subalgebra of locally finite derivations
of $\AA$ such that $\AA$ does not have a proper $\DD$-stable ideal,
we shall present another new class of non-graded infinite
dimensional simple Lie algebras which have some features of Block algebras.
These Lie algebras can also be viewed as further generalizations of Xu's
algebras [15].
\par
Throughout the paper, let $\F$ be a field with
characteristic 0 (in Section 3, we require that $\F$ is algebraically closd).
Denote
the ring of integers by $\Z$ and $\{0,1,2,...\}$ by $\N$.
Up to equivalence of constructions, we can construct the new Lie algebras
as follows. Let $\G$ be an additive subgroup of $\F^4$.
An element in $\G$ is written as $\a=(\a_1,\a_2,\a_3,\a_4)$.
Let $\vf_p:\G\rar\F$ be the natural projection $\vf_p:\a\mapsto\a_p$ for
$p=1,2,3,4$.
We shall assume that $\G$ satisfies the following conditions:
$$
\si\!=\!(1,{\ssc\!}0,{\ssc\!}1,{\ssc\!}0)\in\G,\,
\d\!=\!(0,{\ssc\!}0,{\ssc\!}\d_3,{\ssc\!}0)\in\G, \mbox{ and
}\kn_{\vf_4}\!\not\subset\!\kn_{\vf_2}\mbox{ if } \vf_2\!\ne\!0,
\eqno(1.1)$$%
where $\d_3\!\in\!\F\bs\{0\}$ is a fixed number. Pick
$$
J_p=\{0\}\mbox{ or }\N\mbox{ for }p=1,2,3,4
\mbox{ such that }
J_q\ne\{0\}\mbox{ if }\vf_q=0\ \for\ q=2,4.
\eqno(1.2)$$
Set $J=J_1\times J_2\times J_3\times J_4$.
An element in $J$ will be written as $\ii=(i_1,i_2,i_3,i_4)$.
Denote by $\AA=\AA(\G,J)$ the semigroup algebra of $\G\times J$
with basis $\{x^{\a,\ii}\,|\,(\a,\ii)\in\G\times J\}$ and
the algebraic operation:
$$
x^{\a,\ii} x^{\b,\jj}=x^{\a+\b,\ii+\jj}\ \for\
(\a,\ii),(\b,\jj)\in\G\times J.
\eqno(1.3)$$
Denote the identity element $x^{0,0}$ simply by 1.
Denote
$$
a_{[1]}=(a,0,0,0),\ a_{[2]}=(0,a,0,0),\
a_{[3]}=(0,0,a,0),\ a_{[4]}=(0,0,0,a)\ \for\ a\in\F.
\eqno(1.4)$$
Define derivations $\ptl_p$ of $\AA$ by
$$
\ptl_p(x^{\a,\ii})=\a_p x^{\a,\ii}+i_p x^{\a,\ii-1_{[p]}}\ \for\
(\a,\ii)\in\G\times J\mbox{ and }p=1,2,3,4,
\eqno(1.5)$$
where we always treat $x^{\a,\ii}=0$ if $\ii\notin J$.
\par
We define the following algebraic operation $[\cdot,\cdot]$ on $\AA$
$$
[u,v]\!=\!x^{\si,0}
(\ptl_1\lsb u\rsb\ptl_2\lsb v\rsb
\!-\!\ptl_1\lsb v\rsb\ptl_2\lsb u\rsb)
\!+\!(x^{\d,0}\ptl_3\lsb u\rsb\!+\!u)\ptl_4\lsb v\rsb
\!-\!\ptl_4\lsb u\rsb(x^{\d,0}\ptl_3\lsb v\rsb\!+\!v),
\eqno(1.6)$$
for $u,v\in\AA$.
It can be verified that $(\AA,[\cdot,\cdot])$ forms a Lie algebra. Denote by
$\BB=\BB(\G,J,\d)=[\AA,\AA]$ the derived subalgebras of $(\AA,[\cdot,\cdot])$.
\par
We shall prove the simplicity of $\BB$ and present some interesting
examples in Section 2, then determine the
isomorphism classes and the structure space of the Lie algebras $\BB$ in
Section 3.
\par\
\vs{-15pt}\par
\cl{\bf\S 2. \ Simplicity theorem and examples}
\par
By (1.6), we have
$$
\matrix{
[x^{\a,\ii},x^{\b,\jj}]=\!\!\!\!\!&
(\a_1\b_2-\b_1\a_2)x^{\a+\b+\si,\ii+\jj}
+(\a_1j_2-\b_1i_2)x^{\a+\b+\si,\ii+\jj-1_{[2]}}
\vs{4pt}\hfill\cr&
+(i_1\b_2-j_1\a_2)x^{\a+\b+\si,\ii+\jj-1_{[1]}}
+(i_1j_2-j_1i_2)x^{\a+\b+\si,\ii+\jj-1_{[1]}-1_{[2]}}
\vs{4pt}\hfill\cr&
+(\a_3\b_4-\b_3\a_4)x^{\a+\b+\d,\ii+\jj}
+(\a_3j_4-\b_3i_4)x^{\a+\b+\d,\ii+\jj-1_{[4]}}
\vs{4pt}\hfill\cr&
+(i_3\b_4-j_3\a_4)x^{\a+\b+\d,\ii+\jj-1_{[3]}}
+(i_3j_4-j_3i_4)x^{\a+\b+\d,\ii+\jj-1_{[3]}-1_{[4]}}
\vs{4pt}\hfill\cr&
+(\b_4-\a_4)x^{\a+\b,\ii+\jj}+(j_4-i_4)x^{\a+\b,\ii+\jj-1_{[4]}},
\hfill\cr
}
\eqno(2.1)$$
for $(\a,\ii),(\b,\jj)\in\G\times J$. In particular,
$$
\matrix{
[1,x^{\b,\jj}]\!\!\!\!\hfill&=\!\!\!\!\hfill&
\b_4x^{\b,\jj}+j_4x^{\b,\jj-1_{[4]}}
\mbox{ for }(\b,\jj)\in\G\times J,
\vs{4pt}\hfill\cr
[x^{-\si},x^{\b,\jj}]\!\!\!\!\hfill&=\!\!\!\!\hfill&
-\b_2x^{\b,\jj}-j_2x^{\b,\jj-1_{[2]}}
\mbox{ if }\b_4=j_4=0,
\hfill\cr}
\eqno\matrix{(2.2)\vs{4pt}\cr(2.3)\cr}$$
where in general, we denote $x^\a=x^{\a,0},\,t^{\ii}=x^{0,\ii},\
t_p=t^{1_{[p]}},\,p=1,2,3,4$ (then $x^{\a,\ii}=x^\a t^\ii$ and
$t^\ii=\prod_{p=1}^4t_p^{i_p}$).
For $\l\in\pi_4(\G),m\in\N$, define
$$
\matrix{\dis
\BB^{(m)}_\l=\{u\in\BB\,|\,(\ad_1-\l)^{m+1}(u)=0\}
=\BB\bigcap{\rm span}\{x^{\a,\ii}\,|\,\a_4=\l,i_4\le m\},
\vs{4pt}\hfill\cr\dis
\BB^{(m)}=\bigoplus_{\l\in\pi_4(\G)}\BB^{(m)}_\l,\ \
\BB_\l=\bigcup_{m=0}^\infty\BB^{(m)}_\l,\ \
\mbox{then }\BB=\bigoplus_{\l\in\pi_4(\G)}\BB_\l.
\hfill\cr}
\eqno(2.4)$$
\hs{3ex}
{\bf Theorem 2.1}. {\it The Lie algebra $(\BB,[\cdot,\cdot])$ is simple.
Furthermore, $\BB=\AA$ if $J_p\ne\{0\}$ for some $p\in\{1,2,4\}$, and
if $J_1\!=\!J_2\!=\!J_4\!=\!\{0\}$,
$\BB$ is spanned by the elements:
$$
\{x^{\a,\ii},x^{\b,\ii},
\g_3x^{\g+\d,\ii}\!+\!i_3x^{\g+\d,\ii-1_{[3]}}\!+\!2x^{\g,\ii}\,|\,
(\a_2,\a_4)\!\ne\!0\!=\!(\b_2,\b_4)\!=\!(\g_2,\g_4),
\b_1\!\ne\!1\!=\!\g_1\}.
\eqno(2.5)$$
}\hs{3ex}
{\it Proof}. If $J_4\ne\{0\}$, (2.2) and induction
on $i_4$ show that $x^{\b,\ii}\in[\AA,\AA],\forall
(\b,\ii)\in\G\times J$, so $\BB=\AA$. If $J_4=\{0\}\ne J_2$,
(2.2) shows that
$x^{\b,\ii}\in\BB$ if $\b_4\ne0$, (2.3) and induction on $i_2$ show that
$x^{\b,\ii}\in\BB$ if $\b_4=0$. Assume that $J_2=J_4=\{0\}$. Then
$\pi_2\ne0\ne\pi_4$ by (1.2), and
(2.2), (2.3) show that $x^{\b,\ii}\in\BB$ if $(\b_2,\b_4)\ne0$.
If $(\b_2,\b_4)=0$, for any $\a\in\G$ with $(\a_2,\a_4)\ne0$, we have
$$
[x^{-\a},x^{\a+\b,\ii}]=
\a_2(\b_1x^{\b+\si,\ii}+i_1x^{\b+\si,\ii-1_{[1]}})
+\a_4(\b_3x^{\b+\d,\ii}+i_3x^{\b+\d,\ii-1_{[3]}}+2x^{\b,\ii}).
\eqno(2.6)$$
By the last condition of (1.1), we can choose $\a$ with
$\a_2\!\ne\!0\!=\!\a_4$, thus
the first term of the right-hand side is in $\BB$,
from this we see that $x^{\b,\ii}\!\in\!\BB$ if
$J_1\!\ne\!\{0\}$, i.e., $\BB\!=\!\AA$.
Assume that $J_1\!=\!\{0\}$.
Using (1.2) and the last condition of (1.1),
we see that both terms of the right-hand side of (2.6) are in $\BB$,
from this we see that $x^{\b,\ii}\!\in\!\BB$ if $\b_1\!\ne\!1$ and thus
(2.5) is contained in $\BB$. Conversely, suppose $J_1=J_2=J_4=\{0\}$.
Then the only possible 4 nonzero terms in the right-hand side of (2.1) are
the 1st, 5th, 7th, 9th terms.
If $(\a_2+\b_2,\a_4+\b_4)\ne0$ or
$\a_1+\b_1\ne0,1$, all these terms are in (2.5).
Assume that $(\a_2+\b_2,\a_4+\b_4)=0$ and
$\a_1+\b_1=0,1$.
If $\a_1+\b_1=0$, then the 1st term is zero and the other 3 terms are in
(2.5); if $\a_1+\b_1=1$, then the 1st term and the sum of other 3 term are
in (2.5). Thus
we have (2.5).~This proves the second statement of the theorem.
In particular in any case, we have
$$
x^{i\si},x^{j\d}\in\BB\ \for\ i,j\in\Z,\, i\ne1.
\eqno(2.7)$$
\hs{3ex}
Suppose ${\cal I}\ne\{0\}$ is an ideal of $\BB$. We prove
${\cal I}=\BB$. First we want to prove $1\in{\cal I}$.
Since any ideal is invariant under $\ad_1$, by (2.2),
${\cal I}=\oplus_{\l\in\pi_4(\G)}{\cal I}_\l$, where
${\cal I}_\l={\cal I}\cap\BB_\l$, and
${\cal I}_\l\ne\{0\}\Rla {\cal I}_\l^{(0)}={\cal I}
\cap\BB_\l^{(0)}\ne\{0\}$. Thus ${\cal I}_\l^{(0)}\ne\{0\}$ for some $\l$.
Pick $u\in{\cal I}_\l^{(0)}\bs\{0\}$ and we can write
$$
u=\!\sum_{(\a,\ii)\in M_0}\!c_{\a,\ii}x^{\a,\ii},\mbox{ where }
M_0=\{(\a,\ii)\!\in\!\G\times J\,|\,
\a_4\!=\!i_4\!=\!0,c_{\a,\ii}\!\in\!\F\bs\{0\}\}\mbox{ is finite}.
\eqno(2.8)$$
If $\l\ne0$, by taking $(\b,\jj)\in M_0$ and replacing $u$ by $[x^{-\b},u]$
(it is easy to see that it is nonzero), we can suppose $\l=0$.
Applying $\ad_{x^{-\si}}$ to $u$ (note that
$x^{-\si}\in\BB$ by (2.7)), using (2.3), we see that there exists an
element, again denoted by $u$, such that $\a_2=\eta,i_2=0$ if
$(\a,\ii)\in M_0$, for some $\eta\in\pi_2(\G)$.
\par
\ul{{\it Case 1}: $\pi_2\ne0$.}
If necessary by replacing $u$ by $[x^\b,u]\in{\cal I}\bs\{0\}$
for some $\b\in\G$ with $\b_2\ne0=\b_4$ (cf.~the last condition of (1.1)),
we can suppose $\eta\ne0$.
Denote
$\G(u)=\{\a_1\in\pi_1(\G)\,|\,\exists(\a,\ii)\in M_0\}$, and denote by
$|\G(u)|$ the size of $\G(u)$.
Pick $(\g,\kk)\in M_0$.
Denote $J_\g(u)=\{i_1\in J_1\,|\,\exists(\g,\ii)\in M_0\}$.
Let $v=[x^{\g,\kk},u]$. Note that by (2.1),
$$
[x^{\g,\kk},x^{\a,\ii}]=
\eta((\g_1-\a_1)x^{\a+\g+\si,\ii+\kk}+(k_1-i_1)x^{\a+\g+\si,\ii+\kk-1_{[1]}})
\ \for\ (\a,\ii)\in M_0.
\eqno(2.9)$$
This shows that
$$
\matrix{
\G(v)\!=\!\{\a_1\!+\!\g_1\!+\!1\,|\,\a_1\!\in\!\G(u)\}
\mbox{ if }|J_\g(u)|\!>\!1
\mbox{ or }\{\a_1\!+\!\g_1\!+\!1\,|\,\g_1\!\ne\!\a_1\!\in\!\G(u)\}
\mbox{ if }|J_\g(u)|\!=\!1,
\vs{4pt}\hfill\cr
J_{2\g+\si}(v)=\{i+k-1\,|\,k\ne i\in J_\g(u)\},\hfill\cr}
$$
so $|\G(v)|=|\G(u)|$ if $|J_\g(u)|>1$ or
$|\G(v)|=|\G(u)|-1$ if $|J_\g(u)|=1$, and
$|J_{2\g+\si}(v)|=|J_\g(u)|-1$. Thus by replacing
$u$ by $v$ and repeating the above arguments, we can find an element
$u$ such that $|\G(u)|=1$ and $|J_\g(u)|=1$ for some $\g$, that is,
$M_0=\{(\a,\ii)\in M_0\,|\,\a_1=\g_1,i_1=j\}$ for some $j\in J_1$.
If $j\ne0$, replacing $u$ by $v=[x^\g,u]$ and repeating,
noting that now $[x^\g,x^{\a,\ii}]=-\eta jx^{\a+\g+\si,\ii-1_{[1]}}$
for $(\a,\ii)\in M_0$
(cf.~(2.9)),
we can find an element $u\in{\cal I}$ such that in (2.8),
$M_0=\{(\a,\ii)\in M_0\,|\,(\a_1,\a_2,\a_4)=(\g_1,\eta,0),(i_1,i_2,i_4)=0\}$.
Now taking $(\b,\jj)\in M_0$ and replacing $u$ by $[x^{-\b-\si},u]$, we see
that $u$ has the form (2.8) with $M_0=\{(\a,\ii)\in M_0\,|\,
(\a_1,\a_2,\a_4)=(i_1,i_2,i_4)=0\}$.
\par
\ul{{\it Subcase (i)}: $\pi_4\ne0$.}
By taking $\b\in\G$ with $\b_4\ne0$ and
replacing $u$ by $[x^\b,u]$, we can suppose $M_0\!=\!\{(\a,\ii)\in M_0\,|\,
(\a_1,\a_2,\a_4)\!=\!(\b_1,\b_2,\b_4),(i_1,i_2,i_4)\!=\!0\}$.
Pick $(\g,\kk)\in M_0$ and we have
$$
[x^{\g,\kk},x^{\a,\ii}]=\b_4((\g_3-\a_3)x^{\g+\a+\d,\ii+\kk}+
(k_3-i_3)x^{\g+\a+\d,\ii+\kk-1_{[3]}})
\ \for\ (\a,\ii)\in M_0.
\eqno(2.10)$$
Thus using the same arguments as in (2.9), we can find an element, again
denoted by $u$ such that $M_0$ is a singleton $\{(\a,0)\}$ for some
$\a$ with $\a_4\ne0$, i.e., $u$ is a scalar of $x^\a$.
Thus $1=(2\a_4)^{-1}[x^{-\a},x^\a]\in{\cal I}$.
\par
\ul{{\it Subcase (ii)}: $\pi_4=0$.}
By replacing $u$ by $[t^{2_{[4]}},u]$, we can suppose
$M_0=\{(\a,\ii)\in M_0\,|\,(\a_1,\a_2,\a_4)=0,(i_1,i_2,i_4)=(0,0,1)\}$.
Note that
$$
[x^{\g,\kk},x^{\a,\ii}]\!=\!(\g_3\!-\!\a_3)x^{\g+\a+\d,\ii+\kk-1_{[4]}}\!+\!
(k_3\!-\!i_3)x^{\g+\a+\d,\ii+\kk-1_{[3]}-1_{[4]}}
\ \for\ (\g,\kk),(\a,\ii)\!\in\! M_0.
\eqno(2.11)$$
Thus again using similar arguments as above we obtain $1\in{\cal I}$.
\par
\ul{{\it Case 2}: $\pi_2=0$.}
By replacing $u$ by $[t^{2_{[2]}},u]$, we can suppose
$M_0=\{(\a,\ii)\in M_0\,|\,\a_4=i_4=0,i_2=1\}$. Then similar to
(2.9), we have
$$
[x^{\g,\kk},x^{\a,\ii}]\!=\!(\g_1\!-\!\a_1)x^{\g+\a+\si,\ii+\kk-1_{[2]}}
\!+\!(k_1\!-\!i_1)x^{\g+\a+\si,\ii+\kk-1_{[1]}-1_{[2]}}
\ \for\ (\g,\kk),(\a,\ii)\!\in\! M_0,
\eqno(2.12)$$
and as above, we can deduce $1\in{\cal I}$.
\par
Thus we have $1\in{\cal I}$. Then by (2.2), we have ${\cal I}=\BB$ if
$J_4\ne\{0\}$. Suppose $J_4=\{0\}$. (2.2) shows that $x^{\b,\ii}\in{\cal I}$
if $\b_4\ne0$, and then by letting $\a\in\G$ with $\a_4=0$ and calculating
the coefficient of $k$ in $[x^{\b-(k+1)\si},x^{-\b+k\si+\a,\ii}]\in{\cal I}$,
it gives $\a_2x^{\a,\ii}+i_2x^{\a,\ii-1_{[2]}}\in{\cal I}$. Thus
$x^{\a,\ii}\in{\cal I}$ if $J_2\ne\{0\}$, i.e.,
${\cal I}=\BB$.
Suppose $J_2=\{0\}$. Then $x^{\a,\ii}\in{\cal I}$ if $\a_2\ne0$.
If $\a_2=0$, by taking $\g\in\G$ with $\g_2\ne0=\g_4$ and calculating
the coefficient of $k$ in $[x^{k\g},x^{-k\g+\a-\si,\ii}]\in{\cal I}$, it gives
$(\a_1-1)x^{\a,\ii}+i_1x^{\a,\ii-1_{[1]}}\in{\cal I}$. This shows that
$x^{\a,\ii}\in{\cal I}$ if $J_1\ne\{0\}$, i.e., ${\cal I}=\BB$, and
that $x^{\a,\ii}\in{\cal I}$ if $\a_1\ne1$ and $J_1=\{0\}$. From this and
(2.6), we obtain that elements in (2.5) are all in ${\cal I}$. Thus in any
case ${\cal I}=\BB$, i.e., $\BB$ is simple.
\qed
\par
{\bf Example 2.2}. By Theorem 2.1, we obtain the following examples of
simple Lie algebras $\BB$, where symbol $f_t$ denotes the partial derivative
of a polynomial $f$ with respect to the variable $t$, and where
$m\in\Z\bs\{0\}$ is any fixed number, $a,b,c\in\F\bs\Q$ are any fixed
non-rational numbers in $\F$, and $f,g\in\BB$.
\par\ni
(1) Let $\G=\{(i,0,k,0)\,|\,i,k\in\Z\}, J=\{(0,j,0,l)\,|\,j,l\in\N\}$,
$\d=(0,0,m,0)$. Let $x_1=x^{(1,0,0,0)},x_3=x^{(0,0,1,0)}$,
then $\BB=\F[x^{\pm1}_1,t_2,x^{\pm1}_3,t_4]$ and the Lie bracket is:
$$
[f,g]=x^2_1x_3(f_{x_1}g_{t_2}-f_{t_2}g_{x_1})
+(x_3^{m+1}f_{x_3}+f)g_{t_4}-f_{t_4}(x_3^{m+1}g_{x_3}+g).
\eqno(2.13)$$
(2) Let $\G\!=\!\{(i,j,k,l)\,|\,i,j,k,l\!\in\!\Z\},J\!=\!\{0\},
\d\!=\!(0,0,m,0)$. Let
$\AA\!=\!\F[x^{\pm1}_p\,|\,p\!=\!1,2,3,4]$. Then $\BB\!=\!
\{x_1^{n_1}x_2^{n_2}x_3^{n_3}x_4^{n_4},x_1^{p_1}x_3^{p_3},
x_1(qx_3^{q+m}\!+\!2x_3^q){\sc\,}|{\sc\,}
n_1,...,n_4,p_1,p_3,q\!\in\!\Z,(n_2,n_4)\!\ne\!0,
p_1\!\ne\!1\}$
with the bracket
$$
[f,g]=x^2_1x_2x_3(f_{x_1}g_{x_2}-f_{x_2}g_{x_1})
+x_4((x_3^{m+1}f_{x_3}+f)g_{x_4}-f_{x_4}(x_3^{m+1}g_{x_3}+g)).
\eqno(2.14)$$
(3) Let $\G=\{(i,0,i+ak,0)\,|\,i,k\in\Z\}, J=\{(0,j,0,l)\,|\,j,l\in\N\}$,
$\d=(0,0,a{\ssc\,}m,0)$. Let $x_1=x^{(1,0,1,0)},x_3=x^{(0,0,a,0)}$. Then
$\BB=\F[x_1^{\pm1},t_2,x^{\pm1}_3,t_4]$ with the bracket
$$
[f,g]\!=\!x^2_1(f_{x_1}g_{t_2}\!-\!f_{t_2}g_{x_1})
\!+\!(x_3^m(x_1f_{x_1}\!+\!ax_3f_{x_3})\!+\!f)g_{t_4}
\!-\!f_{t_4}(x_3^m(x_1g_{x_1}\!+\!ax_3g_{x_3})\!+\!g).
\eqno(2.15)$$
(4) Let $\G=\{(i+al,j,i+bj+ck,l)\,|\,i,j,k,l\in\Z\},
J=\{(0,0,0,l)\,|\,l\in\N\}$,
$\d=(0,0,b{\ssc\,}m,0)$. Let $x_1=x^{(1,0,1,0)},x_2=x^{(0,1,b,0)},
x_3=x^{(0,0,c,0)},x_4=x^{(a,0,0,1)}$. Then
$\BB=\F[x_p^{\pm1},t_4\,|\,p=1,2,3,4]$ with the bracket
$$
\matrix{
[f,g]=\!\!\!\!&
x_1x_2((x_1f_{x_1}+ax_4f_{x_4})g_{x_2}-f_{x_2}(x_1g_{x_1}+ax_4g_{x_4}))
\vs{4pt}\hfill\cr&
+(x_3^m(x_1f_{x_1}+bx_2f_{x_2}+cx_3f_{x_3})+f)(x_4g_{x_4}+g_{t_4})
\vs{4pt}\hfill\cr&
-(x_4f_{x_4}+f_{t_4})(x_3^m(x_1f_{x_1}+bx_2g_{x_2}+cx_3g_{x_3})+g).
\hfill\cr}
\eqno(2.16)$$
\par\
\vs{-13pt}\par
\cl{\bf\S 3. \ Isomorphism classes and the structure space}
\par
In this section, we shall assume that $\F$ is an algebraically closed field
with characteristic zero.
To give the structure space, let
$M=\{\vec m=(m_1,m_2,m_3,m_4)\,|\,m_p=0,1,\,p=1,2,3,4\}$ be a set of
16 elements corresponding to all possible choices of $J$.
For $\vec m\in M,\,a\in\F\bs\{0\}$, let $\Omega_{\vec m,a}$
be the set of subgroups $\G$ of $\F^4$, where $\G$ satisfies the following
conditions:
$$
(1,0,1,0),(0,0,a,0)\!\in\!\G;\
\kn_{\vf_4}\!\not\subset\!\kn_{\vf_2}\mbox{ if }
\vf_2\!\ne\!0;\
\vf_p\!\ne\!0\mbox{ if }m_p\!=\!0\ \for\ p\!=\!2,4.
\eqno(3.1)$$
Let $G_{\vec m,a}$ be the group of invertible $4\times4$ matrices of the form
$\pmatrix{^1_{a_1}\, ^0_{a_2}\!\!\!\!&0\cr0\!\!\!\!&
^1_{a_3}\, ^0_{a_{4_{_{}}}}\cr}$
such that $a_1=0$ if $m_1=0\ne m_2$, and $a_3=0$ if $m_3=0\ne m_4$.
Then $G_{\vec m,a}$ acts on $\Omega_{\vec m,a}$ via
$g(\G)=\{\a g^{-1}\,|\,\a\in\G\}\in\Omega_{\vec m,a}$ for $g\in G_{\vec m,a},
\G\in\Omega_{\vec m,a}$.
We define the moduli space ${\cal M}_{\vec m,a}=
\Omega_{\vec m,a}/G_{\vec m,a}$, the set of $G_{\vec m,a}$-orbits in
$\Omega_{\vec m,a}$.
\par
{\bf Theorem 3.1}. {\it The Lie algebras
$\BB=\BB(\G,J,\d)$ and $\BB'=\BB(\G',J',\d')$ are isomorphic $\Rla
(J,\d)=(J',\d')$ and $\exists\,
a_1,a_2,a_3,a_4\in\F$ with $a_2,a_4\ne0$, and $a_1=0$ if $J_1=\{0\}\ne J_2$,
and $a_3=0$ if $J_3=\{0\}\ne J_4$, such that the map
$$\tau:
\b\!=\!(\b_1,\b_2,\b_3,\b_4)\mapsto\b'\!=\!
\b
\pmatrix{^1_{a_1}\, ^0_{a_2}\!\!\!\!&0\cr0\!\!\!\!&^1_{a_3}\, ^0_{a_4}\cr}
\!=\!(\b_1\!+\!a_1\b_2,a_2\b_2,\b_3\!+\!a_3\b_4,a_4\b_4),
\eqno(3.2)$$
is a group isomorphism $\G\cong\G'$. Thus there is a bijection
between the isomorphism classes of the simple Lie algebras $\BB$
and the set:
$$
{\cal M}=\{(\vec m,a,\omega)\,|\,\vec m\in M,a\in\F\bs\{0\},
\omega\in{\cal M}_{\vec m,a}\},
\eqno(3.3)$$
i.e., ${\cal M}$ is the structure space of these simple Lie algebras.
}\par
{\it Proof}.
We shall ALWAYS use the symbol $y$ in place of $x$ for the algebra $\BB'$
and use the same notation with a prime to denote any other element associated
with $\BB'$.\par
``$\Lar$'':
We define two algebra structures (not necessarily associative)
$(\AA,\odot_1)$ and $(\AA,\odot_2)$
by $u\odot_1{\sc}v=x^\si\ptl_1(u)\ptl_2(v)$,
$u\odot_2{\sc}v=(x^\d\ptl_3(u)+u)\ptl_4(v)$,
then $[u,v]=[u,v]_1+[u,v]_2,$ where $[u,v]_p=u\odot_p{\sc}v-v\odot_p{\sc}u$ for
$p\!=\!1,2$.
First define a multiplicative function $\chi:\G\!\rar\!\F\bs\{0\}$ such that
$\chi(\si)\!=\!a_2a_4^{-1},\chi(\d)\!=\!1$. We prove that such a function exists:
suppose $\D$ is the maximal subgroup of $\G$ containing $\{\si,\d\}$
such that such a function
exists for $\D$. If $\D\ne\G$, choose $\b\in\G\bs\D$. If
$\exists\,n\!\in\!\N$ such that $\a\!=\!n\b\!\in\!\D$, then we define
$\chi(\b)
\!=\!\chi(\a)^{1/n}$ (since $\F$ is algebraically closed, we can take any
$n$th root), otherwise we set $\chi(\b)\!=\!1$.
Then $\chi$ extends to a function on the group generated by $\D$
and $\b$, which contradicts the maximality of $\D$.
\par
Define $\th:\AA\rar\AA'$ by
$$
\th(x^{\a,\ii})\!=\!a_4^{-1}\chi(\a)
y^{\a'}(t'_1)^{i_1}(a_1t'_1\!+\!a_2t'_2)^{i_2}
(t'_3)^{i_3}(a_3t'_3+a_4t'_4)^{i_4}, \eqno(3.4)$$ where
$\a'\!=\!\tau(\a)$ for $\a\!\in\!\G$. Note that if
$J_1\!=\!\{0\}\!\ne\!J_2$, by assumption $a_1\!=\!0$ and so $t'_1$
does not appear in $\th(x^{\a,\ii})$; similarly, $t'_3$ does not
occur if $J_3\!=\!\{0\}\!\ne\!J_4$. We want to prove
$$
[\th(x^{\a,\ii}),\th(x^{\b,\jj})]_p=
\th([x^{\a,\ii},x^{\b,\jj}])_p\ \for\ p=1,2.
\eqno(3.5)$$
Say $p=2$ and assume that $i_3+j_3,i_4+j_4>0$ (otherwise, the proof is
easier). We have
$$
x^{\a,\ii}\odot_2 x^{\b,\jj}=
x^{\a+\b,\ii+\jj-1_{[3]}-1_{[4]}}(x^\d(\a_3t_3+i_3)+t_3)(\b_4t_4+j_4).
\eqno(3.6)$$ First we will calculate $\th(x^{\a,\ii}\odot_2
x^{\b,\jj})$. By (3.6) and that $\chi(\d)=1$ and $\tau$ maps $\d$
to $\d$, $\th(x^{\a,\ii}\odot_2 x^{\b,\jj})$ can be written as the
product of the term
$$
a_4^{-1}\chi(\a+\b)y^{\a'+\b'}(t'_1)^{i_1+j_1}
(a_1t'_1+a_2t'_2)^{i_2+j_2}(t'_3)^{i_3+j_3-1}(a_3t'_3+a_4t'_4)^{i_4+j_4-2},
\eqno(3.7)$$
with the term
$$
(a_3t'_3+a_4t'_4)(y^\d(\a_3t'_3+i_3)+t'_3)(\b_4(a_3t'_3+a_4t'_4)+j_4),
\eqno(3.8)$$ (note that $\th(x^{\a+\b,\ii+\jj-1_{[3]}-1_{[4]}})$
is (3.7) multiplied by the first factor of (3.8)). Similarly,
$\th(x^{\a,\ii})\odot_2\th(x^{\b,\jj})$ is (3.7) multiplied by
$$
a_4^{-1}
((y^\d(\a'_3t'_3+i_3)+t'_3)(a_3t'_3+a_4t'_4)+i_4a_3t'_3y^\d)
(\b'_4(a_3t'_3+a_4t'_4)+a_4j_4),
\eqno(3.9)$$
which is arisen from
$(y^\d\ptl'_3+1)(y^{\a'}(t'_3)^{i_3}(a_3t'_3+a_4t'_4)^{i_4})$ and
$\ptl'_4(y^{\b'}(a_3t'_3+a_4t'_4)^{j_4})$.
Denote by $D(\a,\b)$ the difference of (3.8) from (3.9), then
(3.5) is equivalent to $D(\a,\b)-D(\b,\a)=0$, which can be verified
directly by using $\a'_3=\a_3+a_3\a_4,\a'_4=a_4\a_4$.
Thus $\th$ induces an isomorphism $\th:(\BB,[\cdot,\cdot])\cong
(\BB',[\cdot,\cdot])$.
\par
``$\Rar$'': Suppose $\th:\BB\cong\BB'$ is an isomorphism. For a
Lie algebra ${\cal L}$, we denote by ${\cal L}_F$ the set of {\it
ad}-locally finite elements of ${\cal L}$. The isomorphism $\th$
maps $\BB_F$ onto $\BB'_F$. First we determine $\BB_F$.
\par
{\bf Claim 1}. (i) $\BB_F=\F$ if $\vf_4\ne0$; (ii) $\BB_F={\rm
span}\{x^{\a,\ii}\,|\,i_2=i_4=0\}$ if $\pi_2=\vf_4=0$; (iii)
$\BB_F={\rm span}\{x^{-\si},x^{\a,\ii}\,|\,
\a=(\a_3)_{[3]},\ii=(i_3)_{[3]}\}$ (cf.~(1.4)) if
$\pi_2\ne0=\vf_4$; (iv) $\ad_{\sF}$ are locally nilpotent
$\Rla\pi_4=0$.
\par
In each case, the set involved is an abelian Lie subalgebra of $\BB$ and
it is easy to verify that each term in the set is {\it ad}-locally finite,
and that $\ad_{\sF}$ are locally-nilpotent $\Rla\pi_4=0$.
Now suppose $u\in\BB$ is {\it ad}-locally finite.
Write $u=\sum_{(\a,\ii)\in M_0}c_{\a,\ii}x^{\a,\ii},$ where
$M_0\!=\!\{(\a,\ii)\!\in\!\G\!\times\! J\,|\,c_{\a,\ii}\!\ne\!0\}.$
Set $\G(u)\!=\!\{\a\!\in\!\G\,|\,\exists(\a,\ii)\!\in\! M_0\}$. Choose a
total order on $\G$ compatible with the group structure such that
either $\d\!>\!\si\!>\!0$ or $\d\!<\!\si\!<\!0$. This is possible since
$\d,\si$ are $\F$-linear independent. Let $\g$ be the maximal element in
$\G(u)$.
\par
(i) Suppose $\vf_4\!\ne\!0$. If
$\exists\,\a\!\in\!\G(u)$ such that $(\a_3,\a_4)\!\ne\!0$, by reversing the order
of $\G$ if necessary we assume that $(\g_3,\g_4)\!\ne\!0$; otherwise
we assume that $\d\!<\!\si\!<\!0$. Define a total order on $\G\times J:
(\a,\ii)\!<\!(\b,\jj)\Rla\a\!<\!\b$ or $\a\!=\!\b,|\ii|\!<\!|\jj|$
or $(\a,|\ii|)\!=\!(\b,|\jj|)$ but $\exists\,p,\,i_p\!<\!j_p$ and
$i_q\!=\!j_q,q\!<\!p,$
where $|\ii|\!=\!\sum_{p=1}^4i_p$.
Let $(\g,\kk)$ be the maximal element in $M_0$.
We call $c_{\g,\kk}x^{\g,\kk}$ the {\it leading term} of $u$.
We may assume $c_{\g,\kk}=1$. Assume that $(\g,\kk)\ne0$.
If $\d\!<\!\si\!<\!0$, we choose $\eta\!\in\!\G$ such that
$\eta_4\!+\!p\g_4\!\ne\!0,\,\forall\,p\!\ge\!0$,
then the leading term of $\ad_u^q(x^{\eta+\g})$ is
$\prod_{p=0}^{q-1}(\eta_4+p\g_4)x^{\eta+(q+1)\g,q\kk}$,
and ${\rm dim}({\rm span}\{\ad_u^q(x^{\eta+\g})\,|\,q\in\N\})\!=\!
\infty$, so $u$ is not {\it ad}-locally finite.
If $\d\!>\!\si\!>\!0$, then $(\g_3,{\ssc\!}\g_4)\!\ne\!0$, and so we can
choose $\eta\!\in\!\G$ such that
$$
b_p=\g_3(\eta_4+p\g_4)-\g_4(\eta_3+p(\g_3+\d_3))
=\g_3\eta_4-\g_4\eta_3-p\d_3\g_4\ne0, \,\forall\,p\ge0,
\eqno(3.10)$$ then the leading term of $\ad_u^q(x^\eta)$ is
$\prod_{p=0}^{q-1}b_px^{\eta+q(\g+\si),q\kk}$, and so $u$ is not
{\it ad}-locally finite. Thus $\BB_F=\F$.
\par
(ii) Suppose $\vf_2\!=\!\vf_4\!=\!0$. First assume that
$\exists\,(\a,\ii)\!\in\! M_0$
such that $i_4\!\ne\!0$, we define a total order on $\G\times J:
(\a,\ii)\!<\!(\b,\jj)\Rla (i_4,\a,i_1,i_2,i_3)\!<\!(j_4,\b,j_1,j_2,j_3)$
(the later is the lexicographical order).
Let $(\b,\kk)$ be the maximal element in $M_0$. Then $k_4>0$.
If $\exists\,\a\ne0,(\a,\kk)\in M_0$, we choose the order of $\G$
such that $\b>0$; otherwise we choose the order of $\G$ such that $\d<\si<0$.
First suppose $(\b,\kk-1_{[4]})\ne0$.
If $\d<\si<0$, the leading term of $\ad_u^q(t^{(2k_4)_{[4]}})$ is
$\prod_{p=0}^{q-1}(k_4+p(k_4-1))x^{q\b,2\kk+q(\kk-1_{[4]})}$
and $u$ is not {\it ad}-locally finite.
If $\d\!>\!\si\!>\!0$, then $\b\!>\!0$ and we choose $n\in\Z$ such that
$$
b_p=\b_3p(k_4-1)-k_4(n\b_3+p(\b_3+\d_3))=
-nk_4\b_3-p(k_4\d_3+\b_3)\ne0,\,\forall\,p\in\N,
\eqno(3.11)$$
then the leading term of $\ad_u^q(x^{n\b})$ is
$\prod_{p=0}^{q-1}b_px^{n\b+q(\b+\si),q(\kk-1_{[4]})}$
and $u$ is not {\it ad}-locally finite.
If $(\b,\kk)\!=\!(0,1_{[4]})$, if necessary we can reverse the order of $\G$ so
that $\d\!>\!\si\!>\!0$, and still $(0,1_{[4]})$ is the maximal element in
$M_0$,
then the leading term of $\ad_u^q(x^{\si,1_{[4]}})$
is $q!\d_3^qx^{(q+1)\d,1_{[4]}}$, and $u$ is not {\it ad}-locally finite.
Thus $i_4\!=\!0$ for all $(\a,\ii)\!\in\! M_0$. Similarly $i_2\!=\!0$ for
$(\a,\ii)\!\in\! M_0$.
\par
The proof of case (iii) is similar. This proves the claim.
\par
By Claim 1, $\vf_4=0\Rla\vf'_4=0$. We consider the following two cases.
\par
\ul{{\it Case 1}: $\vf_4\ne0$.}
Claim 1 shows that $\th(\F)=\F$. Thus
$$
\th(1)=a\mbox{ for some }a\in\F\bs\{0\}.
\eqno(3.12)$$
So, $\th$ must map $\BB^{(0)}=\oplus_{\l\in\vf_4(\G)}\BB_\l^{(0)}$,
which is the span of eigenvectors of $\ad_1$ (cf.~(2.4)), onto
$\BB'^{(0)}$, and there is a bijection $\tau_4:\l\mapsto\l'$ from $\vf_4(\G)
\rar\vf'_4(\G')$ such that $\th(\BB^{(0)}_\l)=\BB'{}^{(0)}_{\l'}$
and $\th(\BB^{(0)}_0)=\BB'{}^{(0)}_0$.
Since $J_4=\{0\}\Rla\ad_1$ is semi-simple, we have $J_4=J'_4$.
\par
Denote $\CC=\BB^{(0)}_0$. Using (2.3), as in (2.4), for
$\eta\in\pi_2(\G)$, we define
$$
\matrix{\dis
\CC^{(m)}_\eta\!=\!\{u\!\in\!\CC\,|\,
(\ad_{x^{-\si}}\!-\!\eta)^{m+1}(u)\!=\!0\}
\!=\!\BB\bigcap{\rm span}\{x^{\a,\ii}\,|\,
\a_2\!=\!\eta,i_2\!\le\! m,\a_4\!=\!i_4\!=\!0\},
\vs{4pt}\hfill\cr\dis
\CC^{(m)}=\bigoplus_{\eta\in\pi_2(\G)}\CC^{(m)}_\eta,\ \
\CC_\eta=\bigcup_{m=0}^\infty\CC^{(m)}_\eta,\ \
\mbox{then }\CC=\bigoplus_{\eta\in\pi_2(\G)}\CC_\eta.
\hfill\cr}
\eqno(3.13)$$
Note that for the bracket in $\CC$, the last 6 terms of the right-hand side
of (2.1) vanish.
Denote the center of $\CC$ by
$C(\CC)={\rm span}\{x^{\a,\ii}\,|\,\a=(\a_3)_{[3]},\ii=(i_3)_{[3]}\}
\subset\BB$ (cf.~(2.5)).
As in the proof of Claim 1, we have
$$
\matrix{ \CC_F=\F x^{-\si}\!+\!C(\CC) ={\rm
span}\{x^{-\si},x^{\a,\ii}\,|\,\a=(\a_3)_{[3]},\ii=(i_3)_{[3]}\}
\mbox{ if }\vf_2\ne0,\mbox{ or } \vs{4pt}\hfill\cr
\CC_F=\CC_0^{(0)} ={\rm
span}\{x^{\a,\ii}\,|\,(\a_2,\a_4)=(i_2,i_4)=0\} \mbox{ if
}\vf_2=0. \hfill\cr} \eqno(3.14)$$ Thus $\pi_2=0\Rla\CC_F$ has
elements which are not {\it ad}-locally nilpotent on
$\CC\Rla\pi'_2=0$.
\par
\ul{{\it Subcase (i)}: $\pi_2\ne0$.} Then (3.14) shows that
$\pi'_2\ne0,$ and $J_2=J'_2$ since $J_2=\{0\}\Rla \CC_F$ has {\it
ad}-semi-simple elements, and we can suppose
$$
\th(x^{-\si})\!=\!b y^{-\si}\!+\!u_{-\si}\mbox{ for some }b\!\in\!\F\bs\{0\}
\mbox{ and }u_{-\si}\!\in\!C(\CC),
\eqno(3.15)$$
and by (2.3), $\th$ maps $\CC^{(0)}=\oplus_{\eta\in\vf_2(\G)}\CC_\eta^{(0)}$,
onto $\CC'^{(0)}$, and there is a bijection $\tau_2:\eta\mapsto\eta'$ from
$\vf_2(\G)\rar\vf'_2(\G')$ such that $\th(\CC^{(0)}_\eta)=
\CC'{}^{(0)}_{\eta'}$ and $\th(\CC^{(0)}_0)=\CC'{}^{(0)}_0$.
Furthermore, $\th(\CC^{(1)}_0)=\CC'{}^{(1)}_0$. Thus $\th(t_2)=c t'_2+u$ for
some $c\in\F\bs\{0\}$ and $u\in\CC'{}_0^{(0)}.$ Applying $\th$
to $[t_2,x^{-\si}]=1$, by (3.12), (3.15), $c=ab^{-1}$.
Note that $\ptl_{t'_2}:y^{\a',\ii}\!\mapsto\! i_2y^{\a',\ii-1_{[2]}}$ is a
derivation of $\BB'$ if $J'_2\!\ne\!\{0\}$. If $J'_2\!=\!\{0\},$ we set
$\ptl_{t'_2}\!=\!0$.
The isomorphism $\th$ induces a Lie isomorphism of
the derivation algebra of $\BB$ onto the derivation algebra of $\BB'$.
If $J'_2\!\ne\!\{0\}$, we have $(\th^{-1}(\ptl_{t'_2}))(\CC^{(0)})\!=\!
\th^{-1}(\ptl_{t'_2}(\th(\CC^{(0)})))\!=\!\th^{-1}(\ptl_{t'_2}(\CC'^{(0)}))
\!=\!0$ (cf.~(3.13)), and $(\th^{-1}(\ptl_{t'_2}))(t_2)\!=\!
\th^{-1}(\ptl_{t'_2}(\th(t_2)))\!=\!
\th^{-1}(\ptl_{t'_2}(ct'_2\!+\!u))\!=\!\th^{-1}(c)\!=\!b^{-1}$.
Since $\CC$ is generated by $\CC^{(0)}$ and $t_2$, and a derivation is
determined by its actions on generators, we obtain
$\th^{-1}(\ptl_{t'_2})|_{\CC}=b^{-1}\ptl_{t_2}|_{\CC}.$
\par
{\bf Claim 2}. $\th^{-1}(\ptl_{t'_2})=b^{-1}\ptl_{t_2}$.
\par
If $J_2\!=\!\{0\}$, the claim holds trivially. Assume that $J_2\!\ne\!\{0\}$.
So $\BB\!=\!\AA$ by Theorem 2.1. Set
$d\!=\!\th^{-1}(\ptl_{t'_2})\!-\!b^{-1}\ptl_{t_2}$. Then $d(\CC)\!=\!0$.
For $(\b,\jj)\!\in\!\G\!\times\! J$
with $\b_4\!=\!\l\!\ne\!0,j_4\!=\!0$, suppose
$$
d(x^{\b,\jj})=\sum_{(\a,\ii)\in M_{\b,\jj}}b_{\a,\ii}^{(\b,\jj)}
x^{\a+\b,\ii+\jj}\mbox{ for some }b_{\a,\ii}^{(\b,\jj)}\in\F,
\eqno(3.16)$$
where $M_{\b,\jj}=\{(\a,\ii)\in\G\times\Z^4\,|\,\ii+\jj\in J,\,
b_{\a,\ii}^{(\b,\jj)}\ne0\}$ is finite.
Applying $d$ to $[1,x^{\b,\jj}]=\l x^{\b,\jj}$,
using $d(1)=0$, we obtain
$$
\a_4=i_4=0\mbox{ for }(\a,\ii)\in M_{\b,\jj}. \eqno(3.17)$$ Let
$m\!\in\!\Z$ and let $(\g,\kk)\in \G\times J$ be fixed with
$\g_4\!=\!k_4\!=\!0$. Denote
$u_m\!=\!\ad_{x^{-(m+2)\d+\g,\kk}}{\ssc\,}\ad_{x^{m\d}}$. We have
$[x^{m\d},x^{\b,\jj}]=m\d_3\l x^{\b+(m+1)\d,\jj}+\l
x^{\b+m\d,\jj}$, and so
$$
\matrix{
u_m(x^{\b,\jj})=\!\!\!\!\!&
m\d_3\l((\g_1\b_2\!-\!\b_1\g_2)
x^{\b+\g+\si-\d,\jj+\kk}\!+\!(\g_1j_2\!-\!k_2\b_1)x^{\b+\g+\si-\d,\jj+\kk-1_{[2]}}
\vs{4pt}\hfill\cr&
\hs{15pt}+(k_1\b_2\!-\!\g_2j_1)x^{\b+\g+\si-\d,\jj+\kk-1_{[1]}}
\!+\!(k_1j_2\!-\!j_1k_2)x^{\b+\g+\si-\d,\jj+\kk-1_{[1]}-1_{[2]}}
\vs{4pt}\hfill\cr&
\hs{15pt}+(-(m\!+\!2)\d_3\!+\!\g_3)\l x^{\b+\g,\jj+\kk}
\!+\!k_3\l x^{\b+\g,\jj+\kk-1_{[3]}}\!+\!\l x^{\b+\g-\d,\jj+\kk})
\vs{6pt}\hfill\cr&
+\l((\g_1\b_2\!-\!\b_1\g_2)x^{\b+\g+\si-2\d,\jj+\kk}
\!+\!(\g_1j_2\!-\!k_2\b_1)x^{\b+\g+\si-2\d,\jj+\kk-1_{[2]}}
\vs{4pt}\hfill\cr&
\hs{15pt}+(k_1\b_2\!-\!\g_2j_1)x^{\b+\g+\si-2\d,\jj+\kk-1_{[1]}}
\!+\!(k_1j_2\!-\!j_1k_2)x^{\b+\g+\si-2\d,\jj+\kk-1_{[1]}-1_{[2]}}
\vs{4pt}\hfill\cr&
\hs{15pt}+(-(m\!+\!2)\d_3\!+\!\g_3)\l x^{\b+\g-\d,\jj+\kk}
\!+\!k_3\l x^{\b+\g-\d,\jj+\kk-1_{[3]}}
\!+\!\l x^{\b+\g-2\d,\jj+\kk}).
\hfill\cr}
\eqno(3.18)$$
From this, by setting
$v_m=(2\d^2_3\l^2)^{-1}(u_m-2u_{m-1}+u_{m-2})$, we obtain
$v_m(x^{\b,\jj})=x^{\b+\g,\jj+\kk}$ for all
$x^{\b,\jj}\in\BB_\l^{(0)}.$
Applying $d$ to it, noting that $d$ commutes with $v_m$ since $d(\CC)=0$,
by (3.16), (3.17), it gives
$$
b_{\a,\ii}^{(\b,\jj)}=b_{\a,\ii}^{(\b+\g,\jj+\kk)}
\mbox{ for all }(\b,\jj),(\g,\kk)\in\G\times J
\mbox{ with }\b_4=\l,\,j_4=\g_4=k_4=0.
\eqno(3.19)$$
This shows that $b_{\a,\ii}^{(\b,\jj)}$ and $M_{\b,\jj}$ only depend on
$\l$ for all $(\b,\jj)$ with $\b_4\!=\!\l,j_4\!=\!0$. In particular by
setting $\jj\!=\!0$,
it gives $\ii\!\in\! J$ if $(\a,\ii)\!\in\! M_{\b,\jj}$. Set
$w_\l\!=\!\sum_{(\a,\ii)\in M_{\b,\jj}}b_{\a,\ii}^{(\b,\jj)}x^{\a,\ii}\in
\CC$, then (3.16) gives
$$
d(u)=w_\l u\mbox{ for all }u\in\BB_\l^{(0)},
\eqno(3.20)$$
(recall (1.3) for the multiplication on $\AA$). Thus $d^n(u)=w^n_\l u$ for
$n\in\N$. But $d$ is locally nilpotent, we obtain $w_\l=0$,
i.e., $d=0$. If $J_4\ne\{0\}$,
replacing $\jj$ by $\jj+1_{[4]}$ in the discussions above, noting
that all equations hold under mod${\ssc\,}\BB^{(0)}_\l$
(for example, $[1,x^{\b,\jj+1_{[4]}}]\equiv\b_4x^{\b,\jj+1_{[4]}}\,
$(mod${\ssc\,}\BB^{(0)}_\l)$ and (3.18) holds
under mod${\ssc\,}\BB^{(0)}_\l$), we conclude that $d(x^{\b,\jj+1_{[4]}})=0$.
Since $\BB$ is generated by $\BB^{(1)}$, we obtain $d=0$. This proves Claim 2.
\par
Observe that $\ptl^*_2:x^{\b,\jj}\mapsto\b_2x^{\b,\jj}$
for $x^{\b,\jj}\in\BB$, is also a derivation of $\BB$. Similar to the proof
of Claim 2, we obtain
$$
(\th^{-1}(\ptl'{}^*_2))(u)=b^{-1}\ptl^*_2(u)+c_\l u
\mbox{ for }u\in\BB_\l\mbox{ and some }c_\l\in\F.
\eqno(3.21)$$
Denote by $\pi_{24}$ the natural projection $\a\mapsto(\a_2,\a_4)$.
For $(\eta,\l)\in\pi_{24}(\G)$, denote $\BB_{\eta,\l}=\BB\cap{\rm span}
\{x^{\a,\ii}\,|\,\pi_{24}(\a)=(\eta,\l),(i_2,i_4)=0\}$.
Then the nonzero vectors in $\cup_{(\eta,\l)\in\pi_{24}(\G)}\BB_{\eta,\l}$
are the only common eigenvectors of $\ad_1,\ptl_{t_2},\ptl^*_2$. From
(3.12), Claim 2 and (3.21), we see that $\th$ maps common eigenvectors of
$\ad_1,\ptl_{t_2},\ptl^*_2$ to common eigenvectors of
$\ad_1,\ptl_{t'_2},\ptl'{}^*_2$. Thus there is a bijection $\tau_{24}:
\pi_{24}(\G)\rar\pi'_{24}(\G')$ such that $\th(\BB_{\eta,\l})
=\BB'_{\tau_{24}(\eta,\l)}$.
\par
{\bf Claim 3}. If $(\eta,\l)\ne0$, then $\BB_{\eta,\l}$ is a cyclic
module over $\BB_{0,0}=\CC_0^{(0)}$ (cf.~(3.13)), the nonzero scalars of
$x^\a,\,\forall\,\a\in\G$ with $(\a_2,\a_4)=(\eta,\l)$ are the only
generators.
\par
Let $u=\sum_{(\a,\ii)\in M_0}c_{\a,\ii}x^{\a,\ii}\in\BB_{\eta,\l}$,
where $M_0=\{(\a,\ii)\in\G\times J\,|\,(\a_2,\a_4)=(\eta,\l),
i_2=i_4=0, c_{\a,\ii}\ne0\}$ is finite.
If $|M_0|\!>\!1$ or $M_0$ is a singleton $\{(\a,\ii)\}$ with
$\ii\!\ne\!0$, then
$U\!=\!{\rm span}\{x^{\b,\jj}u\!=\!
\sum_{(\a,\ii)\in M_0}c_{\a,\ii}x^{\a+\b,\ii+\jj}\,
|\,(\b,\jj)\!\in\!\G\times J,(\b_2,\b_4)\!=\!(j_2,j_4)\!=\!0\}$
is a proper $\CC_0^{(0)}$-submodule of $\BB_{\eta,\l}$, which contains
$\la u\ra$ (the submodule generated by $u$), so
$u$ is not a generator of $\BB_{\eta,\l}$ as a $\CC^{(0)}_0$-module.
Suppose $M_0=\{(\a,0)\}$ is a singleton with $(\a_2,\a_4)=(\eta,\l)$.
Assume that $\l\ne0$.
By (2.7), $x^{-k\d}\in\CC_0^{(0)}$ for $k\in\N$, we have
$[x^{-k\d},x^\a]=\l(k\d_3x^{\a+(1-k)\d}+x^{\a-k\d})\in\la u\ra,$
thus $x^{\a-k\d}\in\la u\ra,$ and
$$
v_k=[(\g_3+k\d_3)x^{\g+(k+1)\d,\ii}\!+\!i_3x^{\g+(k+1)\d,\ii-1_{[3]}}
\!+\!2x^{\g+k\d,\ii},
x^{\a-(k+2)\d}]\in\la u\ra,
\eqno(3.22)$$
(cf.~(2.5)) for $k\in\N,(\g_2,\g_4)=(i_2,i_4)=0$. Note that $v_k$ is
a polynomial on $k$ with coefficients in $\BB_\l^{(0)}$, computing the
coefficient of $k^2$ shows that $x^{\a+\g,\ii}\in\la u\ra,\,\forall\,
(\g,\vec i)$ with $(\g_2,\g_4)=(i_2,i_4)=0$, but $\BB_{\eta,\l}$ is spanned
by such elements, thus $\BB_{\eta,\l}=\la u\ra$. If $\l=0\ne\eta$,
we obtain the same result by replacing $\d$ by $\si$ in the arguments above.
\par
By Claim 3 and  that $\th(\CC_0^{(0)})=
\CC'{}_0^{(0)}$ (cf.~the statement after (3.15)),
there is a bijection $\tau:\G\bs\kn_{\vf_{24}}
\rar\G'\bs\kn_{\vf'_{24}}$ such that
$$
\th(x^\a)=c_\a y^{\a'}
\mbox{ for }\a\in\G\bs\kn_{\vf_{24}}
\mbox{ and some }c_\a\in\F\bs\{0\},
\eqno(3.23)$$
where when there is no confusion, we denote $\tau(\a)$ by $\a'$
(but we can not denote $\tau(\d)$ by $\d'$ since $\d'$ is the corresponding
element of $\d$ in $\BB'$). Using this and (3.12), we have
$\a_4c_\a y^{\a'}=\th([1,x^\a])=[a,c_\a y^{\a'}]=ac_\a\a'_4  y^{\a'},$
i.e.,
$$\a_4=a\a'_4,
\eqno(3.24)$$
for all $\a\!\in\!\G\bs\kn_{\vf_{24}}$.
For $\a,\b,\a\!+\!\b\!\in\!\G\bs\kn_{\vf_{24}}$,
applying $\th$ to $[x^\a,x^\b]$, by (3.24), we obtain
\par\ni\hs{5pt}$
\matrix{
(\a_1\b_2\!-\!\b_1\a_2)c_{\a+\b+\si}y^{\tau(\a+\b+\si)}
\!\!+\!(\a_3\b_4\!-\!\b_3\a_4)c_{\a+\b+\d}y^{\tau(\a+\b+\d)}
\!\!+\!(\b_4\!-\!\a_4)c_{\a+\b}y^{\tau(\a+\b)}
\vs{4pt}\hfill\cr
=c_\a c_\b((\a'_1\b'_2\!-\!\b'_1\a'_2)y^{\a'+\b'+\si}
\!+\!a^{-1}(\a'_3\b_4\!-\!\b'_3\a_4)y^{\a'+\b'+\d'}
\!+\!a^{-1}(\b_4\!-\!\a_4)y^{\a'+\b'}),
\hfill\cr}
$\hfill(3.25)\par\ni
we obtain that $y^{\tau(\a+\b)}$ must be one of $y^{\a'+\b'+\si},
y^{\a'+\b'+\d'},y^{\a'+\b'},$ i.e.,
$$
\tau(\a+\b)=\a'+\b'+k_{\a,\b}\si+m_{\a,\b}\d',\mbox{ where }
k_{\a,\b},m_{\a,\b}\in\{0,1\}\mbox{ if }\b_4\ne\a_4.
\eqno(3.26)$$
Applying $\th$ to $[x^{-\a},x^{\a}]\!=\!2\a_4$ with $\a_4\!\ne\!0$, denoting
$\b'\!=\!\tau(-\a)$, using (3.12), (3.24), we obtain
$$
c_{-\a}c_\a(\b'_1\a'_2-\a'_1\b'_2)y^{\b'+\a'+\si}
+(\b'_3+\a'_3)a^{-1}\a_4y^{\b'+\a'+\d'}
+2a^{-1}\a_4y^{\b'+\a'})=2a\a_4.
\eqno(3.27)$$
This shows that
$\tau(-\a)=-\a'.$
From this and (3.26),
$(-\a')\!+\!(-\b')\!+\!k_{-\a,-\b}\si\!+\!m_{-\a,-\b}\d'\!=
\!\tau((-\a)\!+\!(-\b))\!=\!
\tau(-(\a\!+\!\b))\!=\!-\tau(\a\!+\!\b)\!
=\!-(\a'\!+\!\b'\!+\!k_{\a,\b}\si\!+\!m_{\a,\b}\d')$, i.e.,
$k_{\a,\b}=m_{\a,\b}=0$. Hence $\tau(\a+\b)\!=\!\tau(\a)\!+\!\tau(\b)$
and $\tau$
can be uniquely extended to a group isomorphism $\tau:\G\!\cong\!\G'$, so
(3.24) holds for all $\a\!\in\!\G$.
By comparing the coefficient of $y^{\a'+\b'}$ in (3.25),
we obtain $c_\a c_\b\!=\!c_{\a+\b}a$ and so $\chi:\a\!\mapsto\! a^{-1}c_\a$ can be
uniquely extended to a multiplicative function: $\G\!\rar\!\F\bs\{0\}$, and the
first terms of both sides of (3.25) then show that $\tau(\si)\!=\!\si$,
similarly, $\tau(\d)\!=\!\d'$.
Fix $\g\!\in\!\kn_{\pi_4}\bs\kn_{\pi_2}$ (cf.~(1.1)). By (3.24) $\g'_4\!=\!0$,
thus by (3.23)
$2\g_2\th(x^{-\si})\!=\!\th([x^{\g-\si},x^{-\g-\si}])
\!=\!c_{\g-\si}c_{-\g-\si}[y^{\g'-\si},y^{-\g'-\si}]
\!=\!2c^2_{-\si}\g'_2y^{-\si}$. Comparing this with (3.15), we obtain
$b\!=\!c^2_{-\si}\g'_2\g_2$ and $u_{-\si}\!=\!0$ in (3.15).
Using this, for any $\a\!\in\!\G$ with $\a_2\!\ne\!0$, applying $\th$ to
$[x^{-\si},x^\a]\!=\!
-\a_2x^\a\!-\!\a_4x^{\a-\si+\d}\!+\!\a_4x^{\a-\si}$,
by (3.23), (3.24), we have $-\a_2c_\a\!=\!-bc_\a\a'_2$, i.e.,
$$
\a_2=b\a'_2,
\eqno(3.28)$$
which must hold for all $\a\in\G$ since $\tau$ is an isomorphism.
Now fix $\b\in\G$ with $\b_2\ne0\ne\b_4$, for $\a\in\G$ with
$\a,\a+\b\in\G\bs\kn_{\pi_{24}}$, using $c_{\a+\b+\si}=c_\a c_\b
c_\si a^{-2}$, then (3.24), (3.28) and the first two terms
of both sides of (3.25) show that
$$
\a'_1=c_\si a^{-2}b\a_1+(-\b_1c_\si a^{-2}b+\b'_1)\b^{-1}_2\a_2,\
\a'_3=c_\d a^{-1}\a_3+(-\b_3c_\d a^{-1}+\b'_3)\b^{-1}_4\a_4,
\eqno(3.29)$$
which must also hold for all $\a\in\G$ since $\tau$ is an isomorphism.
Since $\tau(\si)=\si,$ (3.29) shows that the coefficients of $\a_1,\a_3$
must be 1, i.e., $c_\si a^{-2}b=c_\d a^{-1}=1$. Denote by $a_1,a_3$
the coefficients of $\a_2,\a_4$ respectively, and set $a_2=b^{-1},
a_4=a^{-1}$, then (3.29), (3.24), (3.28) give (3.2) as required.
\par
We already proved $(J_2,J_4)=(J_2',J_4')$. Assume that $J_1\ne\{0\}=J'_1$.
Then we can find $t_1=x^{0,1_{[1]}}\in\CC_0^{(0)}$
(cf. (3.13)) with $[t_1,x^\g]=\g_1 x^{\g+\si}$,
where $\g\in\kn_{\pi_4}\bs\kn_{\pi_2}$ is as above,
but we can not find $u'\in\CC'{}_0^{(0)}$ with
$[u',y^{\g'}]=y^{\g'+\si}$ because for any
$y^{\a',\ii}\in\CC'{}_0^{(0)}$,
$[y^{\a',\ii},y^{\g'}]=(\a'_1\g'_2-\a_2'\g'_1)y^{\g'+\a'+\si}$
can not produce a nonzero term $y^{\g'+\si}$.
This is a contradiction. Thus $J_1=J'_1$.
Similarly $J_3=J'_3$. This proves $J=J'$.
Next we prove $a_1=0$ if $J_1=\{0\}\ne J_2$ as follows.
Consider $\th(x^{\g,1_{[2]}})$. As in the statement after
(3.15), $\th(\CC_\eta^{(1)})=\CC'{}_{\eta'}^{(1)}$,
and since $\ptl_{t_2}(x^{\g,1_{[2]}})=x^\g$, using Claim 2, we can write
$\th(x^{\g,1_{[2]}})=u'+a'y^{\g',1_{[2]}}$ for some
$u'\in\CC'{}^{(0)}_{\eta'}$ and $a'=c_\g b^{-1}$.
Then $\g_1c_{2\g+\si}y^{2\g'+\si}=\th([x^\g,x^{\g,1_{[2]}}])
=c_{\g}([y^{\g'},u']+a'\g'_1 y^{2\g'+\si})$,
but there does not exists $u'\in\CC'{}^{(0)}_{\eta'}$ such that
$[y^{\g'},u']=y^{2\g'+\si}$.
Thus we must have $\g_1c_{2\g+\si}=c_\g a'\g'_1=c_\g^2b^{-1}\g'_1$, i.e.,
$\g_1=\g'_1$. Since $\g'_1=\g_1+a_1\g_2$, it gives
$a_1=0$. Similarly, if $J_3=\{0\}\ne J_4$, then
$a_3=0$. This proves the theorem in this subcase.
\par
\ul{{\it Subcase (ii)}: $\pi_2=0$.}
First a remark: $a_1$ in (3.2) does not affect the
action of $\tau$, so we can always take $a_1\!=\!0$. (3.14) shows that
$\pi'_2\!=\!0$ and $J_2\!=\!J'_2\!=\!\N$ by (1.2).
So $\BB\!=\!\AA,\BB'\!=\!\AA'$ by Theorem 2.1.
(3.14) also shows that $\th(\CC_0^{(0)})\!=\!\CC'{}_0^{(0)}$.
Note that $\CC_0^{(1)}$ is the maximal subspace $U$ of $\CC$
such that $[\CC_0^{(0)},U]\!\subset\!\CC^{(0)}_0$ (cf.~(3.13)), thus we have
$\th(\CC_0^{(1)})\!=\!\CC'{}_0^{(1)}$.
Note that $\CC_0^{(1)}$ is a Lie subalgebra of $\CC$ such that
$\CC_0^{(0)}$ is an abelian ideal of $\CC_0^{(1)}$.
Take the quotient Lie algebra $\ol{\CC}\!=\!\CC_0^{(1)}/\CC_0^{(0)}$.
We use an overbar to denote elements in $\ol\CC$,
then
$$
[\ol x^{\g,\kk+1_{[2]}},\ol x^{\a,\ii+1_{[2]}}]\!=\!(\g_1\!-\!\a_1)
\ol x^{\g+\a+\si,\ii+\kk+1_{[2]}}
\!+\!(k_1\!-\!i_1)\ol x^{\g+\a+\si,\ii+\kk-1_{[1]}+1_{[2]}},
\eqno(3.30)$$
for $(\g,\kk),(\a,\ii)\in\G\times J$ with
$\g_4\!=\!\a_4\!=\!k_2\!=\!k_4\!=\!i_2\!=\!i_4\!=\!0$. As in the proof of
Claim 1, we obtain that the set of {\it ad}-locally
finite elements in $\ol\CC$ is $\F \ol x^{-\si,1_{[2]}}$, and that
$\ol x^{-\si,1_{[2]}}$ is {\it ad}-semi-simple $\Rla J_1=\{0\}$.
Since $\th(\ol\CC)=\ol\CC'$, we obtain $J_1=J'_1$ and we can
write
$$\th(x^{-\si,1_{[2]}})\equiv
a_0y^{-\si,1_{[2]}}\,({\rm mod\,}\CC'{}^{(0)}_0)
\mbox{ for some }a_0\in\F\bs\{0\}.
\eqno(3.31)$$
For $\l\in\pi_1(\G)$, let $\CC^{(0)}_{0,\l}={\rm span}
\{x^{\a,\ii}\in\CC_0^{(0)}\,|\,\a_1=\l,i_1=0\}$, and
$\ol{\CC}_\l={\rm span}\{\ol x^{\a,\ii+1_{[2]}}\in\ol\CC\,|\,
\a_1=\l,i_1=0\}$.
Then $(\cup_{\l\in\pi_1(\G)}\CC^{(0)}_{0,\l})\bs\{0\}$
(resp.~$(\cup_{\l\in\pi_1(\G)}\ol\CC_\l)\bs\{0\}$)
is the set of eigenvectors of $\ol x^{-\si,1_{[2]}}$ in
$\CC_0^{(0)}$ (resp.~$\ol\CC$), we see that there is
a bijection $\tau_p:\pi_1(\G)\rar\pi'_1(\G'),\,p=1,2$ such that
$\th(\CC^{(0)}_{0,\l})=\CC'{}^{(0)}_{0,\tau_1(\l)}$ and
$\th(\ol\CC_\l)=\ol\CC'_{\tau_2(\l)}$.
For $0\ne\l\in\pi_1(\G)$, as in the proof
of Claim 3, we obtain that the nonzero scalars of
$x^\a$ (resp.~$\ol x^{\a,1_{[2]}}$), $\forall\,a\in\G$
with $\a_1=\l,\a_4=0$, are the only generators of
$\CC^{(0)}_{0,\l}$ (resp.~$\ol\CC_\l$) as a $\ol\CC_0$-module.
Thus, there is a bijection $\tau_p:\kn_{\pi_4}\bs\kn_{\pi_1}\rar
\kn_{\pi'_4}\bs\kn_{\pi'_1},p=3,4$ such that
$$
\th(x^\a)\!=\!c_\a y^{\tau_3(\a)},
\ \ \ 
\th(x^{\a,1_{[2]}})\!\equiv\!
c'_\a y^{\tau_4(\a),1_{[2]}}\,({\rm mod\,}\CC'{}^{(0)}_0),
\eqno(3.32)$$
for $\a\!\in\!\kn_{\pi_4}\bs\kn_{\pi_1}$ and some
$c_\a,c'_\a\!\in\!\F\bs\{0\}.$
Applying $\th$ to
$\a_1x^{\a+\b+\si}=[x^\a,x^{\b,1_{[2]}}]$, and to
$(\a_1-\b_1)x^{\a+\b+\si,1_{[2]}}=[x^{\a,1_{[2]}},x^{\b,1_{[2]}}]$,
we obtain
$\tau_3(\a+\b+\si)=\tau_3(\a)+\tau_4(\b)+\si$,
$\a_1c_{\a+\b+\si}=(\tau_3(\a))_1c_\a c'_\b$, and
$\tau_4(\a+\b+\si)=\tau_4(\a)+\tau_4(\b)+\si$,
$(\a_1-\b_1)c'_{\a+\b+\si}=((\tau_4(\a))_1-(\tau_4(\b))_1)c'_\a c'_\b$.
From this we obtain that $\tau_3=\tau_4$ can be uniquely
extended to an isomorphism $\tau_3:\kn_{\pi_4}\!\cong\!\kn_{\pi'_4}$
such that $\tau_3(\si)\!=\!\si$ and by setting $\a\!=\!-\si$, we
obtain $c_\b\!=\!c_{-\si}c'_\b$, and setting $\b\!=\!-\si$ gives $c'_{-\si}\!=\!1$,
i.e., $a_0\!=\!1$ by (3.31). Now as in the arguments before Claim 2, we have
$\th^{-1}(\ptl_{t'_2})|_{\CC_0^{(1)}}\!=\!b^{-1}\ptl_{t_2}|_{\CC_0^{(1)}}$
for some $b$
(a remark: since in this subcase, $\CC$ is not generated by $\CC^{(1)}$,
we can only obtain
$\th^{-1}(\ptl_{t'_2})|_{\CC_0^{(1)}}\!=\!b^{-1}\ptl_{t_2}|_{\CC_0^{(1)}}
\sc\,$).
Then exactly similar to the proof of Claim 2, by letting
$d\!=\!\th^{-1}(\ptl_{t'_2})\!-\!b^{-1}\ptl_{t_2}$, we obtain
$d|_{\BB^{(0)}}\!=\!0$. Since $\BB$ is generated by
$\CC_0^{(1)}$ and $\BB^{(0)}$, we obtain $d=0$.
By considering the common eigenvectors of $\ad_1,\ptl_{t_2}$, we see that
the rest of the proof in this
subcase is the same as that of Subcase (i)
\par
\ul{{\it Case 2}: $\pi_4=0$.}
In this case $a_3$ in (3.2) does not affect the
action of $\tau$, so we can always take $a_3=0$.
By Claim 1, $\pi'_4=0$ and so $J_4=J'_4=\N$ by (1.2).
\par
\ul{{\it Subcase (i)}: $\pi_2=0$} (and so we can take $a_1=0$ in
(3.2)). Denote $\CC_0=\BB_F={\rm span}\{x^{\a,\ii}
\,|\,i_2=i_4=0\},\CC_1=\{u\in\BB\,|\,[u,\CC_0] \subset\CC_0\}={\rm
span}\{x^{\a,\ii}\,|\,i_2+i_4\le1\}, \ol{\CC}\!=\!\CC_1/\CC_0$.
Then $\th(\CC_p)\!=\!\CC'_p,p\!=\!0,1$ and
$\th(\ol{\CC})\!=\!\ol{\CC}'$. Denote by $\ol x^{\a,\ii}$ the
corresponding element of $x^{\a,\ii}$ in $\ol{\CC}$. Since $\CC_0$
is an abelian Lie ideal of $\CC_1$, we can regard $\CC_0$ as a
$\ol{\CC}$-module.
\par
{\bf Claim 4}. The nonzero scalars of $\ol x^{-\si,1_{[2]}}$ are
the only elements in $\ol{\CC}$ acting locally finitely on $\CC_0$.
The adjoint action of $\ol x^{-\si,1_{[2]}}$ on $\CC_0$ is a semi-simple
operator $\Rla J_1=\{0\}$.
\par
Note that for $\a,\b\in\G,\ii,\jj\in J$ with $i_2=i_4=j_2=j_4=0$, we have
$$
\matrix{
[\ol x^{\a,\ii+1_{[2]}},x^{\b,\jj}]
=-\b_1x^{\a+\b+\si,\ii+\jj}-j_1x^{\a+\b+\si,\ii+\jj-1_{[1]}},
\vs{4pt}\hfill\cr
[\ol x^{\a,\ii+1_{[4]}},x^{\b,\jj}]
=-\b_3x^{\a+\b+\d,\ii+\jj}-j_3x^{\a+\b+\d,\ii+\jj-1_{[3]}}-x^{\a+\b,\ii+\jj}.
\hfill\cr}
\eqno(3.33)$$
From this, the proof of the claim is exactly analogous to that of Claim 1.
\par
By Claim 4, $\exists\,a_0\!\ne\!0$ such that
$\th(\ol x^{-\si,1_{[2]}})\!=\!a_0\ol y^{-\si,1_{[2]}}$,
and $J_1\!=\!J'_1$. Let
$\CC_{0,0}\!=\!\{u\!\in\!\CC_0\,|\,[\ol x^{-\si,1_{[2]}},u]\!=\!0\}$,
$\ol{\DD}\!=\!\{u\!\in\!\ol{\CC}\,|\,[\CC_{0,0},u]\!=\!0\}$.
Then by (3.33), $\CC_{0,0}\!=\!{\rm span}\{x^{\b,\ii}\in\CC_0\,|\,
\b_1\!=\!i_1\!=\!i_2\!=\!i_4\!=\!0\}$ and
$\ol{\DD}\!=\!{\rm span}\{\ol x^{\a,\ii+1_{[2]}}\in\ol{\CC}\,|\,
i_1\!=\!i_2\!=\!i_4\!=\!0\}$.
For $\l\!\in\!\pi_1(\G)$, denote $\ol{\DD}_\l
\!=\!{\rm span}\{\ol x^{\a,\ii+1_{[2]}}
\!\in\!\ol{\DD}\,|\,\a_1\!=\!\l,
i_1\!=\!i_2\!=\!i_4\!=\!0\}$. For $(\a,\ii),(\g,\kk)
\!\in\!\G\times J$ with $(i_2,i_4)\!=\!(k_2,k_4)\!=\!0$, we have formula
(3.30),
and we see that $(\cup_{\l\in\pi_1(\G)}\ol{\DD}_\l)\bs\{0\}$ is the set of
eigenvectors of $\ol x^{-\si,1_{[2]}}$ in $\ol{\DD}$. Thus there is
a bijection $\tau_1:\pi_1(\G)\!\rar\!\pi'_1(\G')$ such that
$\th(\ol{\DD}_\l)\!=\!\ol{\DD}'_{\tau_1(\l)}$ for $\l\!\in\!\pi_1(\G)$ and
$\tau_1(0)\!=\!0$.
Using (3.30), as in the proof of Claim 3, we have
\par
{\bf Claim 5}. If $\l\ne0$, then $\ol{\DD}_\l$ is a cyclic module over
$\ol{\DD}_0$,
the nonzero scalars of $\ol x^{\a,1_{[2]}},
\,\forall\,\a\in\G$ with $\a_1=\l$ are the only
generators.
\par
By Claim 5, there is a
bijection $\tau:\G\bs\kn_{\pi_1}\rar\G'\bs\kn_{\pi'_1}$ such that
$$
\th(\ol x^{\a,1_{[2]}})=c_\a\ol y^{\tau(\a),1_{[2]}}
\mbox{ for }\a\in\G\bs\kn_{\pi_1}\mbox{ and some }c_\a\in\F\bs\{0\}.
\eqno(3.34)$$
In particular, $\th$ maps $\ol{\DD}^{(0)}$ to $\ol{\DD}'^{(0)},$
where $\ol{\DD}^{(0)}\!=\!{\rm span}\{\ol x^{\a,1_{[2]}}\!\in\!\ol{\DD}\}$.
Let $\ol{\DD}^{(1)}\!=\!\{u\!\in\!\ol{\CC}\,|\,[u,\ol{\DD}^{(0)}]
\!\subset\! \ol{\DD}^{(0)}\}$.
Then $\th(\ol{\DD}^{(1)})=\ol{\DD}'^{(1)}$. By
$$
[\ol x^{\a,\ii+1_{[4]}},\ol x^{\b,1_{[2]}}]=
\a_1\ol x^{\a+\b+\si,\ii+1_{[4]}}
+i_1\ol x^{\a+\b+\si,\ii-1_{[1]}+1_{[4]}}
-\b_3\ol x^{\a+\b+\d,\ii+1_{[2]}}
-\ol x^{\a+\b,\ii+1_{[2]}},
\eqno(3.35)$$
we obtain $\ol{\DD}^{(1)}=
{\rm span}\{\ol x^{\a,1_{[2]}},\ol x^{\b,1_{[4]}}\,|\,\b_1=0\}$.
Using $[\ol x^{\a,1_{[4]}},\ol x^{\b,1_{[4]}}]
=(\a_3-\b_3)\ol x^{\a+\b+\d,1_{[4]}},$
we see that there is a group isomorphism $\tau_3:\kn_{\pi_1}\cong\kn_{\pi'_1}$
such that
$$
\th(\ol x^{\b,1_{[4]}})=e_\b\ol y^{\tau_3(\b),1_{[4]}}+u_\b
\mbox{ for }\b\in\kn_{\pi_1}\mbox{ and some }e_\b\in\F\bs\{0\},u_\b
\in\ol{\DD}^{(0)}.
\eqno(3.36)$$
Similarly to (3.34), set
$\CC^{(\l)}_0\!=\!\{u\in\CC_0\,|\,[\ol x^{-\si,1_{[2]}},u]\!=\!\l u\}
\!=\!{\rm span}\{x^{\a,\ii}\!\in\!\CC_0\,|\,\a_1\!=\!\l,i_1\!=
\!i_2\!=\!i_4\!=\!0\}$ for $\l\!\in\!\pi_1(\G)$, then as above,
there is a bijection $\tau_0:
{\sc\!}\G\bs\kn_{\pi_1}\!\rar\!\G'\bs\kn_{\pi'_1}$ such that
$$
\th(\ol x^\a)=d_\a\ol y^{\tau_0(\a)}
\mbox{ for }\a\in\G\bs\kn_{\pi_1}\mbox{ and some }d_\a\in\F\bs\{0\}.
\eqno(3.37)$$
Setting $\ii\!=\!\jj\!=\!0$ in (3.33),\,(3.30),\,(3.35), applying $\th$ to
them, using
(3.34),\,(3.36),\,(3.37), we deduce that $\G\!=\!\G'$ and $\tau\!=\!\tau_0$ is
the identity map. To complete the proof of this subcase,
it remains to prove $J_3\!=\!J'_3$, which can be done as in the
last paragraph of Case 1(i).
\par
\ul{{\it Subcase (ii)}: $\pi_2\ne0$.} For $\eta\in\pi_2(\G)$,
denote $\BB_\eta={\rm
span}\{x^{\a,\ii}\,|\,\a_2=\eta,i_2=i_4=0\}$. By Claim 1,
$\BB_F={\rm span}\{x^{-\si},x^{\a,\ii}\,|\,
\a=(\a_3)_{[3]},\ii=(i_3)_{[3]}\}$. Thus
$(\cup_{\eta\in\pi_2(\G)}\BB_\eta)\bs\{0\}$ is the common
eigenvectors of $\BB_F$. So there is a bijection
$\tau_2:\pi_2(\G)\rar\pi'_2(\G')$ such that
$\th(\BB_\eta)=\BB'_{\tau_2(\eta)}$. The rest of the proof is
exactly analogous to that of Case 1(i). \qed\par We remark that as
in the proof of [18, Theorem 4.2], the Lie algebra $\BB(\G,J,\d)$
is not isomorphic to: a Kac-Moody algebra, a generalized Kac-Moody
algebra, a finite root system graded Lie algebra, a generalized
Cartan type Lie algebra\,[8,9,12],\,a Lie algebra of type
L\,[6],\,a Xu's Lie algebra\,[7,11,15],\,a Zhao's Lie algebra in
[18],\,or a Lie algebra in\,[10].
\par
We would like to conclude the paper with the following
\par
{\bf Problem 3.2}. Determine the derivation algebra ${\rm
Der}(\BB)$, the automorphism group ${\rm Aut}(\BB)$, and the
second cohomology group $H^2(\BB,\F)$. \small\par\
\vs{-15pt}\par\ni \cl{\bf References} \vs{-20pt}
\baselineskip=13pt\parskip=7pt
\par\ni\hi4ex\ha1
[1] R.~Block, ``On torsion-free abelian groups and Lie algebras,''
  {\it Proc.~Amer.~Math.~Soc.}, {\bf 9} (1958), 613-620.
  \par\ni\hi4ex\ha1
[2] D.~Dokovic, K.~Zhao, ``Derivations, isomorphisms and
  second cohomology of generalized Block algebras,''
  {\it Algebra Colloquium}, {\bf 3} (1996), 245-272.
  \par\ni\hi4ex\ha1
[3] D.~Dokovic, K.~Zhao, ``Some infinite dimensional simple Lie
  algebras related to those of Block,'' {\it J.~Pure and Applied Algebra},
  {\bf 127} (1998), 153-165.
  \par\ni\hi4ex\ha1
[4] D.~Dokovic, K.~Zhao, ``Some simple subalgebras of generalized
  Block algebras,'' {\it  J.~Alg.}, {\bf 192}, 74-101.
  \par\ni\hi4ex\ha1
[5] J.~M.~Osborn, K.~Zhao, ``Infinite-dimensional Lie algebras
  of generalized Block type,''
  {\it Proc. Amer.~Math.~Soc.}, {\bf 127} (1999), 1641-1650.
  \par\ni\hi4ex\ha1
[6] J.~M.~Osborn, K.~Zhao, ``Infinite-dimensional Lie algebras of type L,''
  {\it Comm.~Alg.}, accepted for publication.
  \par\ni\hi4ex\ha1
[7] Y.~Su, ``Derivations and structure of the Lie algebras of Xu type,''
  {\it Manuscripta Math.}, {\bf105} (2001), 483-500.
  \par\ni\hi4ex\ha1
[8] Y.~Su, X.~Xu, ``Structure of divergence-free Lie algebras,''
  {\it J. Alg.}, {\bf243} (2001), 557-595.

  \par\ni\hi4ex\ha1
[9] Y.~Su, X.~Xu, H.~Zhang, ``Derivation-simple algebras
  and the structures of Lie algebras of Witt type,''
  {\it J. Alg.}, {\bf 233} (2000), 642-662.
  \par\ni\hi4ex\ha1
[10] Y.~Su, J.~Zhou, ``A class of nongraded simple Lie algebras,'' submitted.
  \par\ni\hi4ex\ha1
[11] Y.~Su, J.~Zhou, ``Structure of the Lie algebras related to
  those of Block,'' {\it Comm.~Alg.}, {\bf30} (2002), 3205-3226.
  \par\ni\hi4ex\ha1
[12] X.~Xu, ``New generalized simple Lie algebras of Cartan type over a
  field with characteristic 0,'' {\it J. Alg.}, {\bf 224} (2000), 23-58.
  \par\ni\hi4ex\ha1
[13] X.~Xu, ``On simple Novikov algebras and their irreducible modules,''
  {\it J.~Alg.}, {\bf 185} (1996), 905-934.
  \par\ni\hi4ex\ha1
[14] X.~Xu, ``Novikov-Poisson algebras,'' {\it J.~Alg.}, {\bf190} (1997),
  253-279.
  \par\ni\hi4ex\ha1
[15] X.~Xu, ``Generalizations of Block algebras,'' {\it Manuscripta Math.},
  {\bf 100} (1999), 489-518.
  \par\ni\hi4ex\ha1
[16] X.~Xu, ``Quadratic conformal superalgebras,''
  {\it J.~Alg.}, {\bf 224} (2000), 1-38.
  \par\ni\hi4ex\ha1
[17] X.~Xu, ``Nongraded infinite-dimensional simple Lie algebras,''
  Preprint 2000.
  \par\ni\hi4ex\ha1
[18] K.~Zhao, ``A class of infinite dimensional simple Lie algebras,''
  {\it J.~London Math.~Soc. (2)}, {\bf62} (2000), 71-84.
  \par\ni\hi4ex\ha1
[19] K.~Zhao, ``The $q$-Virasoro-like algebra,'' {\it J.~Alg.}, {\bf188}
  (1997), 506-512.
  \par\ni\hi4ex\ha1
[20] H.~Zhang, K.~Zhao, ``Representations of the Virasoro-like Lie algebra
  and its $q$-analog,'' {\it Comm.~Alg.}, {\bf24} (1996), 4361-4372.
\end{document}